\documentclass[12pt]{article}
\usepackage{float}
\usepackage[german,english]{babel}
\usepackage{epsf,array}
\usepackage{latexsym,bbm}
\usepackage{epsfig,appendix}
\usepackage[]{algorithm2e}
\usepackage {amssymb,amsthm,amsfonts}
\usepackage {amsmath,comment,nccmath}
\usepackage{color,enumitem}
\usepackage{setspace}
\usepackage{url,caption}
\usepackage{hyperref,appendix}
\usepackage[caption=false]{subfig}
\usepackage{url}
\usepackage[noadjust]{cite}
\newcommand{\bm}[1]{\boldsymbol{\mathbf{#1}}}

\usepackage{rotating}
\usepackage[longnamesfirst]{natbib}
\bibpunct{(}{)}{;}{´a}{,}{,}
\renewcommand{\cite}{\citep}

\makeatother
\makeatletter
\renewcommand{\@fnsymbol}[1]{\@arabic{#1}}
\makeatother
\title{On second-order statistics of the log-average periodogram for Gaussian processes}

\author{Karolina Klockmann\footnote{Department of Statistics and Operations Research, Universit\"at Wien,
	Oskar-Morgenstern-Platz 1, 1090 Wien, Austria }
\and Tatyana Krivobokova$^{1}$}

\begin{document}
\maketitle
\begin{abstract}
	\baselineskip=15pt \noindent 
	\\
	We present an approximate expression for the covariance of the log-average periodogram for a zero mean stationary Gaussian process. Our findings extend the work of \citet{ephraim2005second} on the covariance of the log-periodogram  by additionally taking averaging over local frequencies into account. Moreover, we provide a simple expression for the non-integer moments of a non-central chi-squared distribution.
\end{abstract}
\baselineskip=20pt


\section{Introduction} \label{sec:1}

The log-periodogram and the log-average periodogram, where  adjacent frequencies are averaged before applying the logarithm, are fundamental spectral estimators. 
While the second-order properties of the  log-periodogram have been derived in \cite{ephraim2005second}, this  paper investigates the covariance of the log-average periodogram estimator.
Let  $y_1, y_2,\dots, y_{k}$ be a real-valued  time series, which is a realization of a stationary Gaussian process with zero mean.  The spectral density function of the process is assumed to be Hölder continuous. 
Since the time series is  real-valued, we use the discrete cosine transform (DCT) to define the periodogram. This provides estimates of the spectral density on a grid twice as fine compared to the discrete Fourier transform, see \cite{klockmann2024efficient}.
Let
\begin{align}\label{speccomp}
Y_{i}= \frac{1}{\sqrt{k-1}}\left[\frac{y_1+(-1)^iy_{k}}{\sqrt{2}} + \sum_{n=1}^{k-2} y_{n+1} \cos\left (\pi n \frac{i-1}{k-1}  \right)\right ]
\end{align} denote the discrete cosine transform of the time series at frequency $0\leq \pi (i-1)/(k-1) \leq\pi$ for $i=2,...,k-1$. For $i=1,k$, the right-hand side of   (\ref{speccomp}) is additionally multiplied with a factor $1/\sqrt{2}$.  The $Y_i$ are also called spectral components.   The periodogram and  log-periodogram  of the time series are given by $Y_i^2$ and $\log(Y_i^2)$, $i=1,...,k$. 
We define the log-average periodogram  by
\begin{align*}
Y^*_{j}=\log\left ( \frac{1}{m}\sum_{i=(j-1)m+1}^{jm} Y_{i}^2 \right),
\end{align*}  
where $j=1,\dots,\lfloor k/m \rfloor$ and $1\leq m\leq k$ is the number of averaged frequencies.  The number of bins $\lfloor k/m \rfloor $ is denoted by $T$. If $k/m$ is not an integer, the $T$-th bin may contain more than $m$ frequencies.
For finite $k$, the $Y_i,\, i=1,...,k,$ are in general correlated.  We are interested in an expression for the covariance of the log-average periodogram
\begin{align*}
\text{Cov} \left ( Y_j^*, Y_{j^\prime}^*\right), \quad j^\prime, j=1,...,T.
\end{align*}   

The log-average periodogram is widely used in signal processing. While arithmetic averaging is a well-known method to reduce the variance \cite{li2019second}, the logarithm's variance-stabilizing effect for Gamma-distributed data ensures that  $\text{var}[\log(Y_i^2)]$ is constant. Compared to the log-periodogram, additional averaging of adjacent frequencies reduces negative skewness  \cite{heidelberger1981spectral}. The log-average periodogram is for example used to estimate the spectral density, the covariance matrix, the memory (Hurst) parameter, or the cepstral coefficients \cite{klockmann2024efficient,heidelberger1981spectral,robinson1995log,moulines1999broadband,velasco2000non,hurvich2002fexp,sandberg2012optimal}. Several numerical studies suggest that estimators based on the log-average periodogram are more accurate, see \cite{klockmann2024efficient,sandberg2012optimal} and \cite{karagiannis2006understanding}.

Estimation and inference for the periodogram and  estimators based on it usually assume independent spectral components. However, this independence holds only asymptotically  (for $k{\to}\infty$) and for time series with a rapidly decaying autocovariance function. Explicitly considering this dependence allows to improve finite-sample performance and inference for periodogram-based estimators. For example,  \cite{ephraim2005second} derived the covariance of the log-periodogram, showing enhanced performance in cepstral coefficient estimation when accounting for the covariance structure. However, to the best of our knowledge, no formula exists for the covariance matrix of the log-periodogram with averaged frequencies.

This paper derives the covariance structure of the log-average periodogram for a stationary Gaussian process with a smooth spectral density function. The obtained covariance expression extends the work of \cite{ephraim2005second} by accounting for the averaging of adjacent frequencies. 
The new formula is of both theoretical interest and practical use. Furthermore, a new expression for the non-integer moments of a non-central chi-squared distribution is derived, generalizing the results of \cite{rujivan2023analytically} to negative exponents.

\section{Main result} \label{sec:2}
We consider a real-valued stationary Gaussian process with zero mean, covariance matrix   $\bm\Sigma{=}[\sigma_{|i-j|}]_{i,j=1}^k {\in} \mathbb{R}^{k\times k}$ and corresponding spectral density function $f:[-\pi,\pi]\to[0,\infty)$.
The spectral components $Y_i,\, i=1,...,k,$ of the time series are defined as in (\ref{speccomp}). In particular, the vector $(Y_1,...,Y_{k})$ is Gaussian distributed with zero mean  and covariance matrix $\bm{\Omega}{=}[\omega_{st}]_{s,t=1}^k{=}\bm{D}\bm{\Sigma}\bm{D}$ where $\bm{D}$ is the DCT I matrix defined by
\begin{align*} 
\bm{D}&= \sqrt{2} (k-1)^{-1/2}  \left \{\cos \left [\pi (i-1)\frac{j-1}{k-1} \right ]\right\}_{i,j=1}^{k}\\
&\quad \text{ divided by $\sqrt{2}$ when $i,j\in\{1,k\}$.}
\end{align*}

Let $\bm{\Omega}$ be partitioned in $(m\times m)$-dimensional submatrices  as follows:
\begin{align*}    
\bm{\Omega} = \begin{bmatrix} \bm{\Omega}_{1,1} & \cdots & \bm{\Omega}_{1,T} \\
\cdots &\cdots &\cdots\\
\bm{\Omega}_{T,1} & \cdots &\bm{\Omega}_{T,T} \end{bmatrix}.
\end{align*} 
Our result is expressed in terms of the $m$-dimensional correlation matrices, which for $j,j^\prime =1,...,T$, are defined as
$$\bar{\bm{\Omega}}_{j,j^\prime}= \bm{\Omega}_{j,j}^{-1}\bm{\Omega}_{j,j^\prime}\bm{\Omega}_{j^\prime,j}\bm{\Omega}_{j^\prime ,j^\prime}^{-1}.$$
If $k/m$ is not an integer, then $\bm{\Omega}_{j,T}$ with $j\neq T$ is not a square matrix and we define $\bar{\bm{\Omega}}_{j,T}= \bar{\bm{\Omega}}_{T,j}^t= \bm{\Omega}_{j,j}^{-1}\bm{\Omega}_{j,T}\bm{\Omega}_{T ,T}^{-1}\bm{\Omega}_{T,j}.$

Our main result is derived under the following assumptions:
\begin{align*}
& f \text{  is at least } \alpha-\text{Hölder continuous for some } \alpha>0, \label{A1} \tag{A1}\\
& 0< \epsilon \leq f \leq M< \infty \text{ for some } 0<\epsilon<M.\label{A2} \tag{A2}
\end{align*} 
Assumption (\ref{A1}) implies that  $\omega_{st}=\mathcal{O}[\log(k)k^{-\alpha}]$ for $s\neq t=1,...,k$ and $\omega_{st}=f[\pi(s-1)/(k-1)]+\mathcal{O}[\log(k)k^{-\alpha}]$ for $s=t$, where the $\mathcal{O}$ term is uniform over $s,t$ and for $k\to\infty$, see Lemma 1 of \cite{klockmann2024efficient}. 
Assumption (\ref{A2}) implies that the covariance matrix $\bm{\Sigma}$ is invertible.
Under the assumptions (\ref{A1}) and (\ref{A2}), our main result is an approximate expression for the covariance of the log-average periodogram given by
\begin{align} \label{rescov}
&\mbox{Cov}(Y_j^*, Y_{j^\prime}^*)=\sum_{n=1}^\infty\frac{ n!\text{tr}(\bar{\bm{\Omega}}_{j,j^\prime}/m)^n}{(m/2)_nn^2}  [1 + o(1)]\nonumber\\
&=\frac{2\text{tr}(\bar{\bm{\Omega}}_{j,j^\prime}/m)}{m} \,_3F_2(1,1,1,2,m/2+1;\text{tr}(\bar{\bm{\Omega}}_{j,j^\prime}/m))\nonumber\\
&\quad \times  [1 + o(1)],
\end{align} where   $(a)_n{:= }\prod_{i=1}^n(a{+}i{-}1)$ is the Pochhammer symbol, $\, _3F_2(a,b,c,d,e;z){=}\sum_{n=0}^\infty [(a)_n(b)_n(c)_n]/[(d)_n (e)_n]\frac{z^n}{n!}$ is the generalized hypergeometric series \cite[eq. 9.14.1]{lebedev1965special}, and the $o$ term  is with respect to $k \to \infty$ and $m$ fixed.

\section{Derivation of the covariance of the log-average periodogram} \label{sec:3}
Without loss of generality, we consider $j=1$ and $j^\prime=2$. 
The desired covariance expression (\ref{rescov}) is obtained from the second order derivative of the moment generating function
\begin{align*} 
M_{12}(\mu,\eta)&= E\left[ e^{\mu \log (Y_1^2+\dots+Y_{m}^2)+\eta \log(Y_{m+1}^2+\dots+Y_{2m}^2)} \right ] \\
&= E\left [(Y_1^2+\dots+Y_{m}^2)^\mu (Y_{m+1}^2+\dots+Y_{2m}^2)^\eta \right ] 
\end{align*} with respect to $\mu$ and $\eta$ at $\mu=\eta=0$. In particular, $M_{12}(\mu,\eta)\geq 0$. 
Let $\bm{Y}_{j}=(Y_{(j-1)m+1},...,Y_{jm})^t$ for $j=1,2$. 

The moment generating function is evaluated as
\begin{align}
&M_{12}(\mu,\eta)=E\{ (Y_{m+1}^2{+\dots+}Y_{2m}^2)^\eta E [  (Y_1^2{+\dots+}Y_{m}^2)^\mu  | \bm{Y}_2 ] \}. \label{MGF}
\end{align} 
Conditioned on $\bm{Y}_2$ holds  $\bm{Y}_1 \mid \bm{Y}_{2} \sim \mathcal{N}(\bm{\kappa}_{12},{\bm{\Lambda}}_{1,2}) $ where $\bm{\kappa}_{12}=\bm{\Omega}_{1,2}\bm{\Omega}_{2,2 }^{-1} \bm{Y}_{2}$ and  $\bm{\Lambda}_{1,2}=\bm{\Omega}_{1,1}-\bm{\Omega}_{1,2}\bm{\Omega}_{2,2}^{-1}\bm{\Omega}_{2,1}.$
Let $\bm{z}=(z_1,...,z_m)^t$ and $\mbox{d}\bm{z}=\mbox{d}z_1\cdots \mbox{d}z_m$. Then, 
\begin{align}
&E [  (Y_1^2{+\dots+}Y_{m}^2)^\mu  |\bm{Y}_{2} ] 
= \int_{\mathbb{R}^m} (\bm{z}^t\bm{z})^\mu \psi_{\bm{\kappa}_{12},\bm{\Lambda}_{12}}(\bm{z}) \mbox{d}\bm{z} \label{exactE1}
\end{align}  
where $$\psi_{\bm{\kappa},\bm{\Lambda}}(\bm{z}){=}(2\pi)^{-m/2}|\bm{\Lambda}|^{-1/2}\exp[-(\bm{z}{-}\bm{\kappa})^t \bm{\Lambda}^{-1}(\bm{z}{-}\bm{\kappa})/2]$$ is the density function of the multivariate normal distribution with mean vector $\bm{\kappa}$ and covariance matrix $\bm{\Lambda}$.
Let $\bm{U}^t\bm{SU}$ be the eigendecomposition of $\bm{\Lambda}_{1,2}$ with $\bm{S}=\text{diag}(s_1,...,s_m)$. Define  $\lambda_{12}=\text{tr}(\bm{\Lambda}_{1,2})/m$ and let $\mathbb{I}$ denote the identity matrix of dimension $m$. We expand the density function $\psi_{\bm{\kappa}_{12},\bm{\Lambda}_{1,2}}$ with respect to the vector of eigenvalues $\bm{s}{=}(s_1,...,s_m)$ at the point $\bm{\lambda}{=}(\lambda_{12},....,\lambda_{12})$.  The first-order Taylor expansion is given by
\begin{align*}
&\psi_{\bm{\kappa}_{12},\bm{\Lambda}_{1,2}}(\bm{z})=\psi_{\bm{\kappa}_{12},\lambda_{12}\mathbb{I}}(\bm{z})\\
&+\frac{\psi_{\bm{\kappa}_{12},\bm{\Lambda}^*_{12}}(\bm{z})}{2}
\sum_{i=1}^m \left \{\frac{[\bm{U}(\bm{z}-\bm{\kappa}_{12})]_i^2}{(s_i^*)^2}-\frac{1}{s_i^*} \right\}(s_i-\lambda_{12}),
\end{align*} where $\bm{\Lambda}^*_{12}=\bm{U}^t\text{diag}(s^*_1,...,s^*_m)\bm{U}$ for some point $\bm{s}^*=(s_1^*,...,s_m^*)$ between $\bm{s}$ and $\bm{\lambda}$, and  the subscript $i$ denotes the $i$-th component. By assumption (\ref{A1}),  the spectral density function $f$ is $\alpha$-Hölder continuous. Let $\bm\Omega_{1:2,1:2}$ be the $(2m \times 2m)$-dimensional leading  principal submatrix of $\bm\Omega$. Then, by Lemma 1 of \cite{klockmann2024efficient} and the Greshgorin theorem it holds for the set of  eigenvalues $ \lambda(\bm\Omega_{1:2,1:2})$ of $\bm\Omega_{1:2,1:2}$  for some constant $C>0$
$$\lambda(\bm \Omega_{1:2,1:2}) \subset \bigcup_{i=1}^{2m} B\{f[\pi(i-1)/(k-1)],C\log(p)mk^{-\alpha}\},$$
where $B\{a,b\}$ is a closed disc centered at $a$ with radius $b$.

Since $\bm{\Lambda}_{1,2}$ is the Schur complement of block $\bm\Omega_{2,2}$ of matrix $\bm{\Omega}_{1:2,1:2}$, it thus follows: $|s_i-\lambda_{12}|=\mathcal{O}[\log(k)mk^{-\alpha}]$ uniformly over $i$ for $k\to \infty$
and $m$ fixed. Assumption (\ref{A2}) implies $\lambda(\bm{\Lambda}_{1,2})\subset [\epsilon,M]$ and  $s_i^*\in  [\epsilon,M]$.
Let $\Gamma(\cdot)$ denote the Gamma function and let
\begin{align*}
&_0F_1(c;z)=\sum_{n=0}^\infty \frac{1}{(c)_n} \frac{1}{n!} z^n, & _1F_1(b,c;z)=\sum_{n=0}^\infty \frac{(b)_n}{(c)_n} \frac{1}{n!} z^n
\end{align*} denote the hypergeometric series for $b,c\in\mathbb{R}$ such that $c\neq0,-1,...$ and convergence region of $|z|<\infty$  \cite[eq. 9.14.1]{lebedev1965special}.

Define $\bar{\kappa}=\bm{\kappa}_{12}^t\bm{\kappa}_{12}.$ Then, it holds  for $\mu > -m/2$
\begin{align}
& \int_{\mathbb{R}^m} (\bm{z}^t\bm{z})^\mu \psi_{\bm{\kappa}_{12},\lambda_{12} \bm{\mathbb{I}}}(\bm{z}) \mbox{d}\bm{z} \nonumber\\
&{=}\frac{\lambda_{12}^\mu e^{ -\bar{\kappa}/(2\lambda_{12})}}{2^{m/2} \Gamma(m/2)}\int_0^\infty x^{\mu+m/2-1}e^{-x/2}\, _0F_1(m/2;\bar{\kappa}x/(4\lambda_{12})) \mbox{d}x \nonumber\\
&{=}  2^\mu\lambda_{12}^{\mu}  e^{ -\bar{\kappa}/(2\lambda_{12})}  \frac{\Gamma(\mu+m/2)}{\Gamma(m/2)} F_{1,1} (m/2+\mu, m/2;\bar{\kappa}/(2\lambda_{12}))\nonumber \\
&{=}  2^\mu\lambda_{12}^{\mu}   \frac{\Gamma(\mu+m/2)}{\Gamma(m/2)} F_{1,1} (-\mu, m/2;-\bar{\kappa}/(2\lambda_{12})),\label{condexp}
\end{align}
where in the second line a substitution is made to obtain the density function of a non-central chi-squared distribution. In the next two lines,  the identities $\int_0^\infty x^{a-1}e^{-bx}\, dx=b^{-a}\Gamma(a)$ for $a,b>0$ and  $e^{z}\, _1F_1(a,c;-z)=\, _1F_1(c-a,c;z)$ for $a,c\in\mathbb{R}$ and $c\neq0,-1,...,$ are used, see \cite[eq. 3.381]{gradshteyn2014table} and \cite[eq. 9.11.2]{lebedev1965special}. 
In particular, Equation (\ref{condexp}) gives a simple formula for the non-integer moments of a scaled non-central chi-squared distribution. For computational purposes this formula is more convenient than the recently obtained formula by \cite[Cor. 2]{rujivan2023analytically}. Additionally, this formula holds for a range of negative exponents. 
The integral with respect to the first-order remainder term of the Taylor can be bounded by
$ {c_1} \log(k)mk^{-\alpha}\int_{\mathbb{R}^m} [(\bm{z}^t\bm{z})^{\mu+1}{+}\bar{\kappa}(\bm{z}^t\bm{z})^{\mu}]\psi_{\bm{\kappa}_{12},c_2 \mathbb{I}}(\bm{z}) \mbox{d}\bm{z}$ 
for some constants $c_1,c_2>0$. 
Let $\bm{A}_{1,2}=\bm{\Omega}_{2,2}^{-1}\bm{\Omega}_{2,1}\bm{\Omega}_{1,2}\bm{\Omega}_{2,2}^{-1}.$  

Then, we obtain for $\mu>-m/2$
\begin{align}
&M_{12}(\mu,\eta)=  2^\mu\lambda_{12}^{\mu}  \frac{\Gamma(\mu+m/2)}{\Gamma(m/2)} \{ 1+ \mathcal{O}[\log(k)mk^{-\alpha}]\} \nonumber\\ 
&{\times}\int_{\mathbb{R}^m} (\bm{z}^t\bm{z})^\eta F_{1,1} (-\mu, m/2;-\bm{z}^t\bm{A}_{1,2}\bm{z}/(2\lambda_{12})) \psi_{0,\bm{\Omega}_{2,2}}(\bm{z}) \mbox{d}\bm{z} ,\label{MGF2} 
\end{align} where the $\mathcal{O}$ term is with respect to $k\to\infty$ and $m$ fixed.
Let $\bm{B}_{1,2}=\bm{\Omega}_{2,2}^{1/2} \bm{A}_{1,2}\bm{\Omega}_{2,2}^{1/2}$. Next, we insert the integral form of the confluent hypergeometric function  \cite[eq. 9.11.1]{lebedev1965special} , which holds for $z,b,c\in\mathbb{R}$  such that $c>b>0$. For the moment, we assume that $0 < -\mu < m/2$. Then, 
\begin{align}
&\int_{\mathbb{R}^m} (\bm{z}^t\bm{z})^\eta F_{1,1} (-\mu, m/2;-\bm{z}^t\bm{A}_{1,2}\bm{z}/(2\lambda_{12})) \psi_{0,\bm{\Omega}_{2,2}}(\bm{z})\,\mbox{d}\bm{z}\nonumber\\
&= \frac{\Gamma(m/2)}{\Gamma(-\mu)\Gamma(m/2+\mu)}\int_{0}^1t^{-\mu-1}(1-t)^{m/2+\mu-1}\nonumber \\
& \quad \times \int_{\mathbb{R}^m} \frac{(\bm{z}^t\bm{\Omega}_{2,2}\bm{z})^\eta}{(2\pi)^{m/2}} \exp\left [-\frac{\bm{z}^t \left (\mathbb{I}+\bm{B}_{1,2}t/\lambda_{12} \right)\bm{z}}{2}\right ] \, \mbox{d}\bm{z}\, \mbox{d}t. \label{E2_Taylor1} 
\end{align}
Let $\bm{U}^t\bm{S U}$ be the eigendecomposition of $\bm{\Omega}_{2,2}$, where $\bm{S}=\text{diag}(s_1,...,s_m)$. 
Define  $\omega_{2}=\text{tr}(\bm{\Omega}_{2,2})/m$. A first order Taylor expansion of the  function $g(\bm{s})=[\bm{z}^t\bm{U}^t\text{diag}(s_1,...,s_m)\bm{U}\bm{z}]^\eta $ at $\bar{\bm{\omega}}=(\omega_2,....,\omega_2)$ gives
\begin{align*}
g(\bm{s})&= \omega_2^\eta (\bm{z}^t\bm{z})^\eta\\
&\quad +\eta[\bm{z}^t\bm{U}^t\text{diag}(s^*_1{,...,}s^*_m)\bm{U}\bm{z}]^{\eta-1}\sum_{i=1}^m (\bm{Uz})_i^2(s_i{-}\omega_{2}),
\end{align*} where  $\bm{s}^*=(s_1^*,...,s_m^*)$ is between $\bm{s}$ and $\bar{\bm{\omega}}$. 
By assumptions (\ref{A1}) and (\ref{A2}) it holds $\lambda(\bm{\Omega}_{2,2})\subset [\epsilon,M]$ and $|s_i-\omega_2|=\mathcal{O}[\log(k)mk^{-\alpha}]$. 
Furthermore, let $\bm{V}^t\bm{R V}$ be the eigendecomposition of $\bm{B}_{1,2}$, where $\bm{R}=\text{diag}(r_1,...,r_m)$. Define  $\gamma_{12}=\text{tr}(\bm{B}_{1,2})/m$.  A first order Taylor expansion of the  function $h(\bm{r})=\exp[-\bm{z}^t\bm{V}^t\text{diag}(r_1,...,r_m)\bm{V}\bm{z}t/(2\lambda_{12})]$ at $\bar{\bm{r}}=(\gamma_{12},...,\gamma_{12})$ gives
\begin{align*}
h(\bm{r})&=\exp\left(-\bm{z}^t\bm{z}\frac{\gamma_{12}}{2\lambda_{12}}t\right)-\frac{t}{2\lambda_{12}}h(\bm{r}^*) \sum_{i=1}^m(\bm{V}\bm{z})_i^2(r_i-\gamma_{12}),
\end{align*}  where  $\bm{r}^*=(r_1^*,...,r_m^*)$ is between $\bm{r}$ and $\bar{\bm{r}}$. 
The eigenvalues  of $\bm{B}_{1,2}$ are non-negative and the largest eigenvalue is bounded by $c[\log(k)mk^{-\alpha}]^2$
as the entries of  $\bm{\Omega}_{2,1}\bm{\Omega}_{1,2}$ are of the order $\mathcal{O}\{m[\log(k)k^{-\alpha}]^2\}$. In particular, it follows: $|r_i-\gamma_{12}|=\mathcal{O}\{m^2[\log(k)k^{-\alpha}]^2\}$.
Thus, for $\eta>-m/2$
\begin{align}
&\int_{\mathbb{R}^m} \frac{(\bm{z}^t\bm{\Omega}_{2,2}\bm{z})^\eta}{(2\pi)^{m/2}} \exp\left [-\frac{\bm{z}^t \left (\mathbb{I}+\bm{B}_{1,2}t/\lambda_{12} \right)\bm{z}}{2}\right ] \, \mbox{d}\bm{z} \nonumber \\
&= \omega_2^{\eta} \int_{\mathbb{R}^m}\frac{(\bm{z}^t\bm{z})^{\eta}}{(2\pi)^{m/2}} e^{-\frac{\bm{z}^t\bm{z}}{2}\left(1+\frac{\gamma_{12}t}{\lambda_{12}}\right)} \, \mbox{d}\bm{z}\{1 {+} \mathcal{O}[\log(k)mk^{-\alpha}]\} \nonumber\\
&=\frac{\Gamma(\eta+m/2)}{\Gamma(m/2)} \frac{2^\eta\omega_2^{\eta} }{(1+\gamma_{12}t/\lambda_{12})^{\eta+m/2}}\{1 {+}  \mathcal{O}[\log(k)mk^{-\alpha}]\},\label{E2_Taylor2}
\end{align} where the integral was solved with spherical coordinates. 
The remainder terms of the Taylor expansions  can be bounded by  
$c_3\log(k)mk^{-\alpha}\int_{\mathbb{R}^m} (\bm{z}^t\bm{z})^{\eta} \psi_{0,c_4\mathbb{I}}(\bm{z})\,\mbox{d}\bm{z}$ for some constant $c_3,c_4>0$. 
Let $$_2F_1(a,b,c;z)=\sum _{n=0}^{\infty }{\frac {(a)_{n}(b)_{n}}{(c)_{n}}}\,{\frac {z^{n}}{n!}}$$ be the hypergeometric series for $a,b,c\in\mathbb{R}$ such that $c\neq0,-1,...$, and which converges absolutely for $|z|<1$ and for $z=1$ if $c-a-b>0$ \cite[eq. 9.1.2]{lebedev1965special}. 
Using the integral version of $_2F_1(a,b,c;z)$ \cite[eq. 9.1.6]{lebedev1965special}, which holds for $z<1$ and $a,b,c\in\mathbb{R}$ such that $c>b>0$, and plugging (\ref{E2_Taylor2}) and (\ref{E2_Taylor1}) into (\ref{MGF2}), gives
\begin{align}
M_{12}(\mu,\eta)= &2^{\mu+\eta}{\lambda}_{12}^\mu \omega_2^\eta  \frac{\Gamma(\mu + m/2)}{\Gamma(m/2)} \frac{\Gamma\left ( \eta+m/2\right)}{\Gamma(m/2)} \nonumber\\
& \times \, _2F_1(\eta+m/2,-\mu,m/2;-\gamma_{12}/\lambda_{12}) \nonumber\\
& \times \{1 + \mathcal{O}[\log(k)mk^{-\alpha}]\}. \label{MGF3}
\end{align} 

Formula (\ref{MGF3}) was proven under the assumption that $m/2>-\mu>0$. Since $_2F_1(a,b,c;z)/\Gamma(c)$ is an entire function of the parameters $a,b,c$, it can be analytically extended for $\mu>-m/2$ \cite[Ch. 9.4]{lebedev1965special}. 
Let $\omega_1=\text{tr}(\bm{\Omega}_{1,1}/m).$ 
Applying the transformation $$\, _2F_1(a,b,c;z)=(1-z)^{-b}\,_2F_1(c-a,b,c;z/(z-1)),$$  \cite[eq. 9.5.2]{lebedev1965special}  with $z/(z-1)=-\gamma_{12}/\lambda_{12}$  and using the Lipschitz continuity of $\,_2F_1(a,b,c;z)$leads to
\begin{align*}
&M_{12}(\mu,\eta) = 2^{\mu+\eta}\omega_1^\mu\omega_2^\eta  \frac{\Gamma(\mu + m/2)}{\Gamma(m/2)} \frac{\Gamma\left ( \eta+m/2\right)}{\Gamma(m/2)}\\
& \times \,  _2F_1(-\eta,-\mu,m/2;\text{tr}(\bar{\bm{\Omega}}_{1,2}/m))\{1+ \mathcal{O}[\log(k)mk^{-\alpha}]\}.
\end{align*} 

Note that  $\text{tr}(\bar{\bm{\Omega}}_{1,2}/m)\leq 1$ since $0 \leq \text{tr}(\bm{\Omega}_{1,1}^{-1}\bm{\Lambda}_{12}/m)$.
The second order derivative of $M_{12}(\mu,\eta)$ at $\mu=\eta=0$ can be obtained by using $\partial(\mu)_n/\partial \mu= -(n-1)!$ at $\mu=0$ for $n>0$. 
Some calculations using the product rule results in
\begin{align*}
\left .\frac{\partial^2 M_{12}(\mu,\eta)}{\partial \eta \partial \mu}\right|_{\substack{\mu=0\\\eta=0}} 
&{=}\left . \sum_{n=1}^\infty \frac{ \frac{\partial}{\partial \eta} ({-}\eta)_n \frac{\partial}{\partial \mu} ({-}\mu)_n   \text{tr}(\bar{\bm{\Omega}}_{1,2}/m)^n}{(m/2)_n n!} \right|_{\substack{\mu=0\\ \eta=0}}\\
&\quad \times \{1+ \mathcal{O}[\log(k)mk^{-\alpha}]\}\\
&{=}\sum_{n=1}^\infty\frac{ n!\text{tr}(\bar{\bm{\Omega}}_{1,2}/m)^n}{(m/2)_nn^2} \{1+ \mathcal{O}[\log(k)mk^{-\alpha}]\},
\end{align*} where $_2F_1(a,0,c;z)=\,_2F_1(0,b,c;z)=1$ has been used. Since all integrals are finite, the moment generating function exists and is infinitely differentiable with respect to its parameters. This implies that the derivative of  the $\mathcal{O}$ term evaluated at $\mu=\eta=0$ exists and goes to zero at least as fast as $\log(k)mk^{-\alpha}$ for $k\to \infty$. 
Finally,
\begin{align*}
&\sum_{n=1}^\infty\frac{ n!\text{tr}(\bar{\bm{\Omega}}_{1,2}/m)^n}{(m/2)_nn^2}{=}\frac{2}{m}\sum_{n=1}^\infty\frac{(1)_{n-1}(1)_{n-1}(1)_{n-1}}{(m/2+1)_{n-1}(2)_{n-1}}\frac{\text{tr}(\bar{\bm{\Omega}}_{1,2}/m)^n}{(n-1)!}\\
&=\frac{2\text{tr}(\bar{\bm{\Omega}}_{1,2}/m)}{m} \,_3F_2(1,1,1,2,m/2+1;\text{tr}(\bar{\bm{\Omega}}_{1,2}/m)),
\end{align*} since for $n\geq1$ holds $(n-1)!=(1)_{n-1}$, $n=(2)_{n-1}/(1)_{n-1}$ and $(m/2)_n=(m/2) (m/2+1)_{n-1}$.
\section{Application}
\begin{figure}[h] 
	\vspace*{-.2cm}
	\captionsetup{farskip=0pt,nearskip=4pt} 	 
	\subfloat[]{ \includegraphics[width=0.5\textwidth]{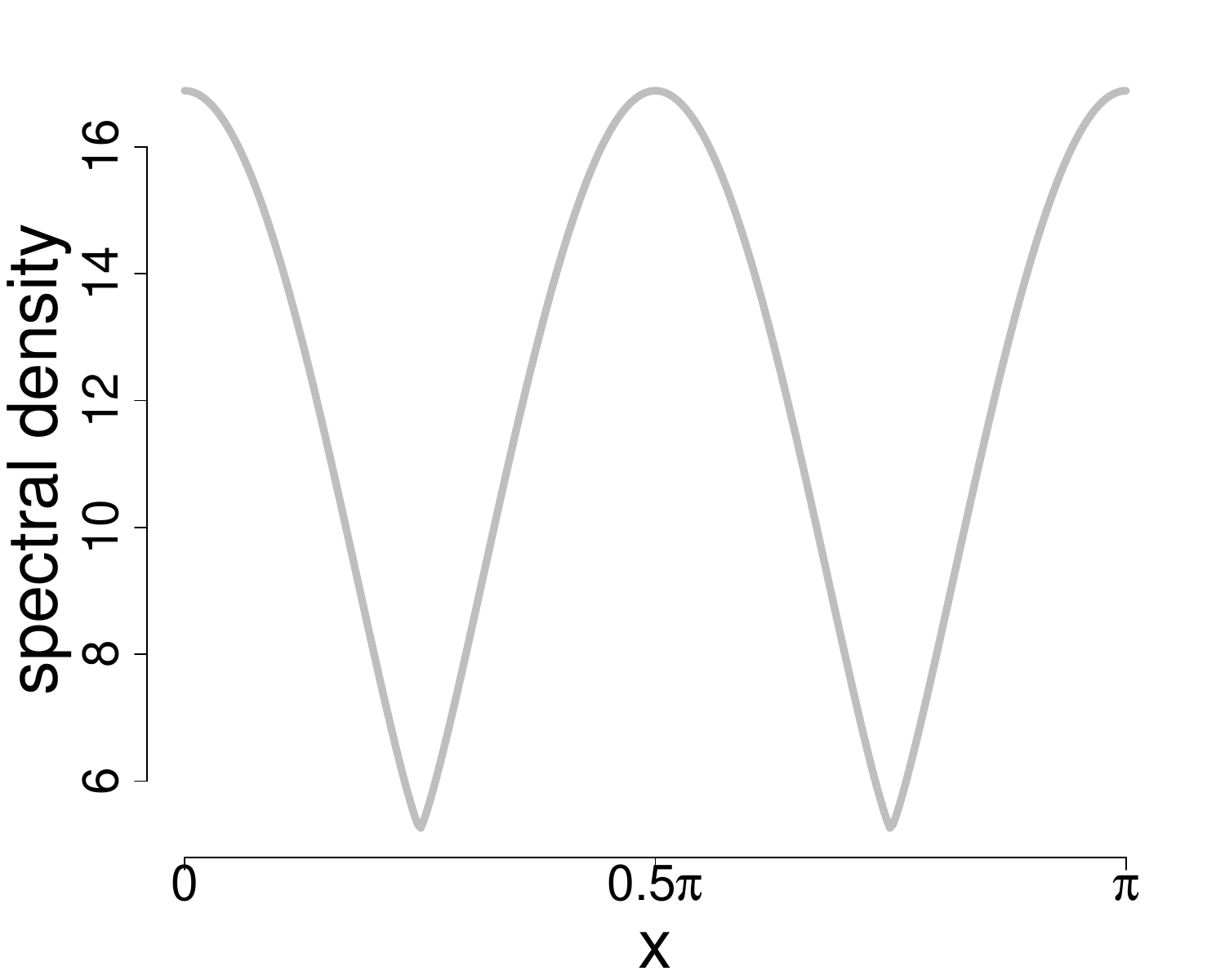}}
	\subfloat[]{\includegraphics[width=0.5\textwidth]{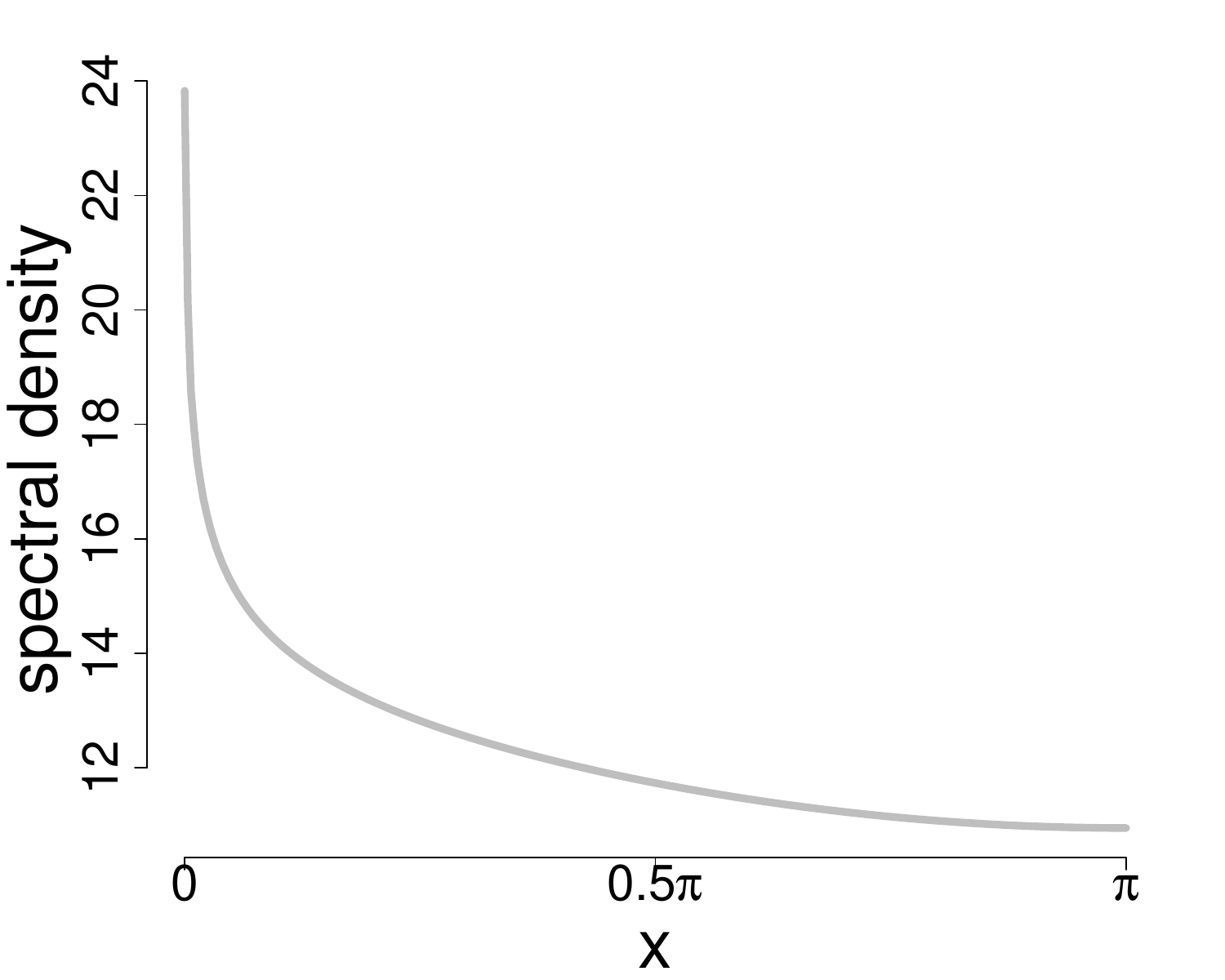}}
	\caption{Spectral density functions for (a) example 1 and (b) example 2. }
	\label{fig_Simulation}
	\vspace*{-.2cm}
\end{figure}

The covariance formula (\ref{rescov}) can be applied in practical situations to improve estimators which are either based on the log-average periodogram, or more general, which require good estimators of the two fundamental quantities, spectral density $f$ and autocovariance function $\sigma$, especially in settings with small time series length $k$. For example, to check stationarity of a Markov chain, the Heidelberger-Welch  test of \cite{heidelberger1981spectral} involves estimating $f(0)$ by fitting a polynomial regression to the first $K$ frequencies of the log-average periodogram. To estimate the Hurst parameter, \cite{kettani2006novel} requires estimates of  $\sigma_0$ and $\sigma_1$. 
The effect of the covariance formula (\ref{rescov}) is investigated for these two examples through a simulation study.
We consider a Gaussian time series of length $p=60$ and with the autocovariance functions $\sigma$ 
\begin{enumerate}
	\item[\textbf{(1)}] such that $f$ is Lipschitz continuous but not differentiable: $f(x)= 3.5^2\{|\sin(x+0.5\pi)|^{1.2}+0.45\},$ 
	\item[\textbf{(2)}] such that $\sigma_i= 3.5^2 (|i+1|^{1.1} - 2 |i|^{1.1} + |i-1|^{1.1})/2.$ 
\end{enumerate}
Figure \ref{fig_Simulation} shows the spectral densities for the two examples.
 A Monte Carlo simulation with $500$ iterations is performed using R (version 4.1.2).  In example \textbf{(1)}, the following methods to estimate $f(0)$ are compared: the  Heidelber-Welch (HW) estimator, the spectral density estimator  VST-DCT of \cite{klockmann2024efficient}, and their corrected versions based on the decorrelated log-average periodogram. 
To compute the covariance matrix (\ref{rescov}), the uncorrected HW and the VST-DCT estimator were employed. Using the true covariance matrix $\Sigma$ gives similar results. The R package \texttt{vstdct} was utilized for the VST-DCT estimator. Cubic splines were selected, with the smoothing parameter chosen with generalized cross-validation. The number of bins was set to  $T=15$, following the method described in \cite{klockmann2024efficient}. 
For the HW estimator, the number of bins is $T=k/2$. 
The fit depends heavily on the polynomial degree $d$ and the number of frequencies $K$ used in the regression. We found that $d=1$ and $K=10$ produced the best results, as shown in Table~\ref{Tab_ex1}, which also includes the error for $f$ and $\Sigma$ w.r.t. $L_\infty$ and spectral norm.  However, as there are no data-driven methods, this outcome is not practically attainable.
In example \textbf{(2)},  the Hurst parameter $H$ of a fractional Gaussian noise process is estimated with  Kettani's method  \cite{kettani2006novel}, which involves estimation of  $\sigma_1$. It is compared to the Hurst estimator which uses $\hat{\sigma}_1$ implied by the VST-DCT estimator of $f$; with and without decorrelating the log-average periodogram, and with the same parameter choices as in the other setting. 
\begin{table}
	\caption{HW and VST-DCT estimator with and without decorrelation for estimating  $f(0)$, $f$ and $\Sigma$ in example \textbf{(1)}.}
	\label{Tab_ex1}
	\small
	\setlength{\tabcolsep}{3pt}
	\vspace*{-.4cm}
	\begin{center}
		\begin{tabular}{ccccc}
			\hline
			& $f(0)$& $\hat{f}(0)$& $\|\hat{f}-f\|_\infty$ & $\|\hat{\bm{\Sigma}}-\bm{\Sigma}\|_2$ \\ 
			\hline
			HW & 16.7984 &  15.8068 & -- & --\\ 
			HW with decor. & 16.7984 & 15.8074 & -- & --\\
			VST-DCT & 16.7984 & 16.5199  & 2.0968  & 1.7937 \\ 
			VST-DCT with decor. & 16.7984 & 16.6093  & 2.0845 & 1.7878 \\  
			\hline
		\end{tabular}
	\end{center}
\end{table}
\begin{table}
	\vspace*{-.4cm}
	\caption{Kettani's and VST-DCT estimator with and without decorrelation for estimating   $\sigma_1$  and $ H$ in example \textbf{(2)}.}
	\small
	\setlength{\tabcolsep}{3pt}
	\vspace*{-.4cm}
	\begin{center}
		\begin{tabular}{ccccc}
			\hline
			&$\sigma_1$& $\hat{\sigma}_1$ &$H$ &$\hat{H}$\\ 
			\hline
			Kettani & 0.8792 &0.5255 & 0.55& 0.5244 \\ 
			VST-DCT  & 0.8792  &0.6928& 0.55& 0.5301 \\ 		
			VST-DCT with decor.  & 0.8792 &0.7196 &0.55 & 0.5312 \\ 
			\hline
		\end{tabular}
		\vspace*{-.6cm}
	\end{center}
\end{table}
The results show that the VST-DCT estimators  perform better in both examples and could be further improved with the new decorrelation formula (\ref{rescov}). As the HW estimator is based on a parametric model, it is expected that the decorrelation formula has a smaller effect. 
All in all, the simulation study demonstrated the superiority of the estimators for $f$ and $\sigma$ based on the decorrelated log-average periodogram, particularly in scenarios with shorter signal lengths. Since these quantities are involved in numerous testing and estimation procedures for time series, our findings are highly relevant in practice.  

Further applications include covariance and precision  matrix estimation \cite{klockmann2024efficient}, cepstrum estimation \cite{sandberg2012optimal}, estimation of the long-run variance \cite{hannan1957variance} or  testing a time series for stationarity \cite{paparoditis2009testing}.
\section{Discussion}
We have derived an approximate expression for the covariance of the log-average periodogram for stationary Gaussian time series with an invertible covariance matrix and a bounded, $\alpha$-Hölder continuous spectral density.  
The main challenge in deriving expression (\ref{rescov}) is computing non-integer moments of a weighted sum of non-central chi-squared random variables in (\ref{exactE1}) and the mixed moments of a weighted sum of chi-squared random variables in (\ref{E2_Taylor1}). While \cite[Cor. 2]{rujivan2023analytically} provided an analytic expression for non-integer moments of a weighted sum of non-central chi-squared random variables, it only holds for positive exponents.  For (\ref{E2_Taylor1}), there is currently no known analytic expression. Using a Taylor expansion, we simplify (\ref{exactE1}) and (\ref{E2_Taylor1}) to an unweighted form, deriving  formula (\ref{condexp}) for non-integer moments of the non-central chi-squared distribution. This formula is simpler, for computational purposes more convenient and also holds for negative exponents.

The  obtained covariance formula (\ref{rescov})  shows that the variance of the log-average periodogram is  $\psi(m/2) [1+o(1)]$  where $\psi$ is the trigamma function, and the covariance of two different components of the log-average periodogram is of the order $\mathcal{O}[\log(k)^2k^{-2\alpha}]$ for $k{\to} \infty$ and $m$ fixed. 
For $m{=}1$ our results match those of \cite{ephraim2005second} up to a factor of $0.5$ due to the use of the discrete cosine transform. The derivation of  (\ref{rescov}) with the discrete Fourier transform results in the same formula without the $0.5$ factor. 
Finally, we conjecture that a similar formula could be derived for Gaussian time series with only positive semidefinite covariance matrices by carefully handling the quadratic forms involving the degenerate Gaussian distribution and using generalized matrix inversions.

\bibliographystyle{apalike}
\bibliography{Literature}

\begin{thebibliography}{}

\bibitem[Ephraim and Roberts, 2005]{ephraim2005second}
Ephraim, Y. and Roberts, W.~J. (2005).
\newblock On second-order statistics of log-periodogram with correlated
  components.
\newblock {\em IEEE Signal Process. Lett.}, 12(9):625--628.

\bibitem[Gradshteyn and Ryzhik, 2014]{gradshteyn2014table}
Gradshteyn, I.~S. and Ryzhik, I.~M. (2014).
\newblock {\em Table of integrals, series, and products}.
\newblock Academic press.

\bibitem[Hannan, 1957]{hannan1957variance}
Hannan, E. (1957).
\newblock The variance of the mean of a stationary process.
\newblock {\em Journal of the Royal Statistical Society Series B: Statistical
  Methodology}, 19(2):282--285.

\bibitem[Heidelberger and Welch, 1981]{heidelberger1981spectral}
Heidelberger, P. and Welch, P.~D. (1981).
\newblock A spectral method for confidence interval generation and run length
  control in simulations.
\newblock {\em Communications of the ACM}, 24(4):233--245.

\bibitem[Hurvich et~al., 2002]{hurvich2002fexp}
Hurvich, C.~M., Moulines, E., and Soulier, P. (2002).
\newblock The fexp estimator for potentially non-stationary linear time series.
\newblock {\em Stochastic processes and their applications}, 97(2):307--340.

\bibitem[Karagiannis et~al., 2006]{karagiannis2006understanding}
Karagiannis, T., Molle, M., and Faloutsos, M. (2006).
\newblock Understanding the limitations of estimation methods for long-range
  dependence.
\newblock {\em University of California}.

\bibitem[Kettani and Gubner, 2006]{kettani2006novel}
Kettani, H. and Gubner, J.~A. (2006).
\newblock A novel approach to the estimation of the long-range dependence
  parameter.
\newblock {\em IEEE Transactions on Circuits and Systems II: Express Briefs},
  53(6):463--467.

\bibitem[Klockmann and Krivobokova, 2024]{klockmann2024efficient}
Klockmann, K. and Krivobokova, T. (2024).
\newblock Efficient nonparametric estimation of toeplitz covariance matrices.
\newblock {\em Biometrika}, page asae002.

\bibitem[Lebedev et~al., 1965]{lebedev1965special}
Lebedev, N.~N., Silverman, R.~A., and Livhtenberg, D. (1965).
\newblock Special functions and their applications.
\newblock {\em Physics Today}, 18(12):70.

\bibitem[Li et~al., 2019]{li2019second}
Li, T., Fan, H., Garc{\'\i}a, J., and Corchado, J.~M. (2019).
\newblock Second-order statistics analysis and comparison between arithmetic
  and geometric average fusion: Application to multi-sensor target tracking.
\newblock {\em Information Fusion}, 51:233--243.

\bibitem[Moulines and Soulier, 1999]{moulines1999broadband}
Moulines, E. and Soulier, P. (1999).
\newblock Broadband log-periodogram regression of time series with long-range
  dependence.
\newblock {\em The Annals of Statistics}, 27(4):1415--1439.

\bibitem[Paparoditis, 2009]{paparoditis2009testing}
Paparoditis, E. (2009).
\newblock Testing temporal constancy of the spectral structure of a time
  series.
\newblock {\em Bernoulli}, 15(4):1190--1221.

\bibitem[Robinson, 1995]{robinson1995log}
Robinson, P.~M. (1995).
\newblock Log-periodogram regression of time series with long range dependence.
\newblock {\em The Annals of Statistics}, pages 1048--1072.

\bibitem[Rujivan et~al., 2023]{rujivan2023analytically}
Rujivan, S., Sutchada, A., Chumpong, K., and Rujeerapaiboon, N. (2023).
\newblock Analytically computing the moments of a conic combination of
  independent noncentral chi-square random variables and its application for
  the extended cox--ingersoll--ross process with time-varying dimension.
\newblock {\em Mathematics}, 11(5):1276.

\bibitem[Sandberg and Hansson-Sandsten, 2012]{sandberg2012optimal}
Sandberg, J. and Hansson-Sandsten, M. (2012).
\newblock Optimal cepstrum smoothing.
\newblock {\em Signal process.}, 92(5):1290--1301.

\bibitem[Velasco, 2000]{velasco2000non}
Velasco, C. (2000).
\newblock Non-gaussian log-periodogram regression.
\newblock {\em Econometric Theory}, 16(1):44--79.

\end{thebibliography}

\end{document}